\documentclass[11pt]{article}
\usepackage{amsmath, upgreek}
\usepackage{amssymb}
\usepackage{amscd}
\usepackage{xspace}
\usepackage{verbatim}
\usepackage{color}
\usepackage{latexsym}
\usepackage{times,url}
\usepackage{amsmath,amssymb}
\usepackage{amsfonts}
\usepackage{graphicx}
\usepackage{subfigure}
\usepackage{indentfirst}
\newcommand{\mysection}[1]{
\section{#1}\setcounter{equation}{0}}
\title{\bf General uniqueness results for large solutions}
\author{{\bf Juli\'an L\'opez-G\'omez\footnote{  Instituto de Matem\'atica Interdisciplinar (IMI), Departamento de Analisis Matem\'atico y Matem\'atica Aplicada,
Universidad Complutense de Madrid, 28040-Madrid, Spain.\newline
Email: Lopez-Gomez@mat.ucm.es}}\qquad
{\bf Luis Maire\footnote{  Departamento de Matem\'atica Aplicada, Ciencia e Ingenier\'ia de los Materiales y Tecnolog\'ia Electr\'onica, Spain. \newline
Email: luis.maire@urjc.es}}
\qquad
 {\bf Laurent V\'eron\footnote{ D\'epartement de Math\'ematiques,  Universit\'e Fran\c{c}ois Rabelais, Tours, France.
 \newline Email: veronl@lmpt.univ-tours.fr }}\\[2mm]
}

\date{}
\begin{document}
 \maketitle


\newcommand{\txt}[1]{\;\text{ #1 }\;}
\newcommand{\tbf}{\textbf}
\newcommand{\tit}{\textit}
\newcommand{\tsc}{\textsc}
\newcommand{\trm}{\textrm}
\newcommand{\mbf}{\mathbf}
\newcommand{\mrm}{\mathrm}
\newcommand{\bsym}{\boldsymbol}
\newcommand{\scs}{\scriptstyle}
\newcommand{\sss}{\scriptscriptstyle}
\newcommand{\txts}{\textstyle}
\newcommand{\dsps}{\displaystyle}
\newcommand{\fnz}{\footnotesize}
\newcommand{\scz}{\scriptsize}
\newcommand{\be}{\begin{equation}}
\newcommand{\bel}[1]{\begin{equation}\label{#1}}
\newcommand{\ee}{\end{equation}}
\newcommand{\eqnl}[2]{\begin{equation}\label{#1}{#2}\end{equation}}
\newcommand{\barr}{\begin{eqnarray}}
\newcommand{\earr}{\end{eqnarray}}
\newcommand{\bars}{\begin{eqnarray*}}
\newcommand{\ears}{\end{eqnarray*}}
\newcommand{\nnu}{\nonumber \\}
\newtheorem{subn}{\name}
\renewcommand{\thesubn}{}
\newcommand{\bsn}[1]{\def\name{#1}\begin{subn}}
\newcommand{\esn}{\end{subn}}
\newtheorem{sub}{\name}[section]
\newcommand{\dn}[1]{\def\name{#1}}   
\newcommand{\bs}{\begin{sub}}
\newcommand{\es}{\end{sub}}
\newcommand{\bsl}[1]{\begin{sub}\label{#1}}
\newcommand{\bth}[1]{\def\name{Theorem}
\begin{sub}\label{t:#1}}
\newcommand{\blemma}[1]{\def\name{Lemma}
\begin{sub}\label{l:#1}}
\newcommand{\bcor}[1]{\def\name{Corollary}
\begin{sub}\label{c:#1}}
\newcommand{\bdef}[1]{\def\name{Definition}
\begin{sub}\label{d:#1}}
\newcommand{\bprop}[1]{\def\name{Proposition}
\begin{sub}\label{p:#1}}
\newcommand{\BBR}{\eqref}
\newcommand{\rth}[1]{Theorem~\ref{t:#1}}
\newcommand{\rlemma}[1]{Lemma~\ref{l:#1}}
\newcommand{\rcor}[1]{Corollary~\ref{c:#1}}
\newcommand{\rdef}[1]{Definition~\ref{d:#1}}
\newcommand{\rprop}[1]{Proposition~\ref{p:#1}}
\newcommand{\BA}{\begin{array}}
\newcommand{\eA}{\end{array}}
\newcommand{\BAN}{\renewcommand{\arraystretch}{1.2}
\setlength{\arraycolsep}{2pt}\begin{array}}
\newcommand{\BAV}[2]{\renewcommand{\arraystretch}{#1}
\setlength{\arraycolsep}{#2}\begin{array}}
\newcommand{\BSA}{\begin{subarray}}
\newcommand{\eSA}{\end{subarray}}
\newcommand{\Remark}{\note{Remark}}
\newcommand{\BAL}{\begin{aligned}}
\newcommand{\eAL}{\end{aligned}}
\newcommand{\BALG}{\begin{alignat}}
\newcommand{\eALG}{\end{alignat}}
\newcommand{\BALGN}{\begin{alignat*}}
\newcommand{\eALGN}{\end{alignat*}}
\newcommand{\note}[1]{\textit{#1.}\hspace{2mm}}
\newcommand{\Proof}{\note{Proof}}
\newcommand{\qeda}{\hspace{10mm}\hfill $\square$}
\newcommand{\qed}{\\
${}$ \hfill $\square$}
\newcommand{\BBRemark}{\note{Remark}}
\newcommand{\modin}{$\,$\\[-4mm] \indent}
\newcommand{\forevery}{\quad \forall}
\newcommand{\set}[1]{\{#1\}}
\newcommand{\setdef}[2]{\{\,#1:\,#2\,\}}
\newcommand{\setm}[2]{\{\,#1\mid #2\,\}}
\newcommand{\mt}{\mapsto}
\newcommand{\lra}{\longrightarrow}
\newcommand{\lla}{\longleftarrow}
\newcommand{\llra}{\longleftrightarrow}
\newcommand{\Lra}{\Longrightarrow}
\newcommand{\Lla}{\Longleftarrow}
\newcommand{\Llra}{\Longleftrightarrow}
\newcommand{\warrow}{\rightharpoonup}
\newcommand{
\paran}[1]{\left (#1 \right )}
\newcommand{\sqbr}[1]{\left [#1 \right ]}
\newcommand{\curlybr}[1]{\left \{#1 \right \}}
\newcommand{\abs}[1]{\left |#1\right |}
\newcommand{\norm}[1]{\left \|#1\right \|}
\newcommand{
\paranb}[1]{\big (#1 \big )}
\newcommand{\lsqbrb}[1]{\big [#1 \big ]}
\newcommand{\lcurlybrb}[1]{\big \{#1 \big \}}
\newcommand{\absb}[1]{\big |#1\big |}
\newcommand{\normb}[1]{\big \|#1\big \|}
\newcommand{
\paranB}[1]{\Big (#1 \Big )}
\newcommand{\absB}[1]{\Big |#1\Big |}
\newcommand{\normB}[1]{\Big \|#1\Big \|}
\newcommand{\produal}[1]{\langle #1 \rangle}

\newcommand{\thkl}{\rule[-.5mm]{.3mm}{3mm}}
\newcommand{\thknorm}[1]{\thkl #1 \thkl\,}
\newcommand{\trinorm}[1]{|\!|\!| #1 |\!|\!|\,}
\newcommand{\bang}[1]{\langle #1 \rangle}
\def\angb<#1>{\langle #1 \rangle}
\newcommand{\vstrut}[1]{\rule{0mm}{#1}}
\newcommand{\rec}[1]{\frac{1}{#1}}
\newcommand{\opname}[1]{\mbox{\rm #1}\,}
\newcommand{\supp}{\opname{supp}}
\newcommand{\dist}{\opname{dist}}
\newcommand{\myfrac}[2]{{\displaystyle \frac{#1}{#2} }}
\newcommand{\myint}[2]{{\displaystyle \int_{#1}^{#2}}}
\newcommand{\mysum}[2]{{\displaystyle \sum_{#1}^{#2}}}
\newcommand {\dint}{{\displaystyle \myint\!\!\myint}}
\newcommand{\q}{\quad}
\newcommand{\qq}{\qquad}
\newcommand{\hsp}[1]{\hspace{#1mm}}
\newcommand{\vsp}[1]{\vspace{#1mm}}
\newcommand{\ity}{\infty}
\newcommand{\prt}{\partial}
\newcommand{\sms}{\setminus}
\newcommand{\ems}{\emptyset}
\newcommand{\ti}{\times}
\newcommand{\pr}{^\prime}
\newcommand{\ppr}{^{\prime\prime}}
\newcommand{\tl}{\tilde}
\newcommand{\sbs}{\subset}
\newcommand{\sbeq}{\subseteq}
\newcommand{\nind}{\noindent}
\newcommand{\ind}{\indent}
\newcommand{\ovl}{\overline}
\newcommand{\unl}{\underline}
\newcommand{\nin}{\not\in}
\newcommand{\pfrac}[2]{\enfrac{(}{)}{}{}{#1}{#2}}

\def\ga{\alpha}     \def\gb{\beta}       \def\gg{\gamma}
\def\gc{\chi}       \def\gd{\delta}      \def\ge{\epsilon}
\def\gth{\theta}                         \def\vge{\varepsilon}
\def\gf{\phi}       \def\vgf{\varphi}    \def\gh{\eta}
\def\gi{\iota}      \def\gk{\kappa}      \def\gl{\lambda}
\def\gm{\mu}        \def\gn{\nu}         \def\gp{\pi}
\def\vgp{\varpi}    \def\gr{\rho}        \def\vgr{\varrho}
\def\gs{\sigma}     \def\vgs{\varsigma}  \def\gt{\tau}
\def\gu{\upsilon}   \def\gv{\vartheta}   \def\gw{\omega}
\def\gx{\xi}        \def\gy{\psi}        \def\gz{\zeta}
\def\Gg{\Gamma}     \def\Gd{\Delta}      \def\Gf{\Phi}
\def\Gth{\Theta}
\def\Gl{\Lambda}    \def\Gs{\Sigma}      \def\Gp{\Pi}
\def\Gw{\Omega}     \def\Gx{\Xi}         \def\Gy{\Psi}

\def\CS{{\mathcal S}}   \def\CM^+{{\mathcal M}}   \def\CN{{\mathcal N}}
\def\CR{{\mathcal R}}   \def\CO{{\mathcal O}}   \def\CP{{\mathcal P}}
\def\CA{{\mathcal A}}   \def\CB{{\mathcal B}}   \def\CC{{\mathcal C}}
\def\CD{{\mathcal D}}   \def\CE{{\mathcal E}}   \def\CF{{\mathcal F}}
\def\CG{{\mathcal G}}   \def\CH{{\mathcal H}}   \def\CI{{\mathcal I}}
\def\CJ{{\mathcal J}}   \def\CK{{\mathcal K}}   \def\CL{{\mathcal L}}
\def\CT{{\mathcal T}}   \def\CU{{\mathcal U}}   \def\CV{{\mathcal V}}
\def\CZ{{\mathcal Z}}   \def\CX{{\mathcal X}}   \def\CY{{\mathcal Y}}
\def\CW{{\mathcal W}} \def\CQ{{\mathcal Q}}
\def\BBA {\mathbb A}   \def\BBb {\mathbb B}    \def\BBC {\mathbb C}
\def\BBD {\mathbb D}   \def\BBE {\mathbb E}    \def\BBF {\mathbb F}
\def\BBG {\mathbb G}   \def\BBH {\mathbb H}    \def\BBI {\mathbb I}
\def\BBJ {\mathbb J}   \def\BBK {\mathbb K}    \def\BBL {\mathbb L}
\def\BBM {\mathbb M}   \def\BBN {\mathbb N}    \def\BBO {\mathbb O}
\def\BBP {\mathbb P}   \def\BBR {\mathbb R}    \def\BBS {\mathbb S}
\def\BBT {\mathbb T}   \def\BBU {\mathbb U}    \def\BBV {\mathbb V}
\def\BBW {\mathbb W}   \def\BBX {\mathbb X}    \def\BBY {\mathbb Y}
\def\BBZ {\mathbb Z}

\def\GTA {\mathfrak A}   \def\GTB {\mathfrak B}    \def\GTC {\mathfrak C}
\def\GTD {\mathfrak D}   \def\GTE {\mathfrak E}    \def\GTF {\mathfrak F}
\def\GTG {\mathfrak G}   \def\GTH {\mathfrak H}    \def\GTI {\mathfrak I}
\def\GTJ {\mathfrak J}   \def\GTK {\mathfrak K}    \def\GTL {\mathfrak L}
\def\GTM {\mathfrak M}   \def\GTN {\mathfrak N}    \def\GTO {\mathfrak O}
\def\GTP {\mathfrak P}   \def\GTR {\mathfrak R}    \def\GTS {\mathfrak S}
\def\GTT {\mathfrak T}   \def\GTU {\mathfrak U}    \def\GTV {\mathfrak V}
\def\GTW {\mathfrak W}   \def\GTX {\mathfrak X}    \def\GTY {\mathfrak Y}
\def\GTZ {\mathfrak Z}   \def\GTQ {\mathfrak Q}

\font\Sym= msam10 
\def\SYM#1{\hbox{\Sym #1}}
\newcommand{\bdw}{\prt\Gw\xspace}

\renewcommand{\thefigure}{\arabic{section}.\arabic{figure}}

\newcommand{\field}[1] {\mathbb{#1}}
\newcommand{\N}{\field{N}}
\newcommand{\Z}{\field{Z}}
\newcommand{\Q}{\field{Q}}
\newcommand{\C}{\field{C}}
\newcommand{\K}{\field{K}}
\newcommand{\fin}{\mbox{$\quad{}_{\Box}$}}
\def\a{\alpha}
\def\b{\beta}
\def\e{\varepsilon}
\def\D{\Delta}
\def\d{\delta}
\def\g{\gamma}
\def\G{\Gamma}
\def\l{\lambda}
\def\m{\mu}
\def\n{\~n}
\def\L{\Lambda}
\def\o{\omega}
\def\O{\Omega}
\def\p{\partial}
\def\r{\rho}
\def\vr{\varrho}
\def\S{\Sigma}
\def\s{\sigma}
\def\t{\theta}
\def\T{\Theta}
\def\v{\varphi}
\def\ov{\overline}
\def\un{\underline}
\def\ua{\uparrow}
\def\da{\downarrow}
\def\ra{\rightarrow}
\def\na{\nearrow}
\def\z{\zeta}

\date{}
\maketitle\medskip

\noindent{\small {\bf Abstract} We give a series of very general sufficient conditions in order to ensure the uniqueness of large solutions
 for $-\Gd u+f(x,u)=0$ in a bounded domain $\Gw$ where $f:\overline\Gw\ti\BBR\mapsto \BBR_+$ is a continuous function, such that $f(x,0)=0$ for $x\in\overline\Gw$, and $f(x,r)>0$ for $x$ in a neighborhood of $\prt\Gw$ and all $r>0$.
\medskip

\vspace{1mm}
\hspace{.05in}
\tableofcontents
\medskip
\mysection{Introduction}
Let  $\Gw \sbs \BBR^N$ be a bounded domain and $f:\overline\Gw\ti\BBR\mapsto \BBR_+$  a continuous function such that $f(x,0)=0$ and $r\mapsto f(x,r)$ is nondecreasing  for $x\in\overline\Gw$, and  $f(x,r)>0$ for $x$ in a neighborhood of $\prt\Gw$ and all $r>0$. This paper deals with the uniqueness question of the solution of the equation
\bel{I.1}\begin{array} {lll}
-\Gd u+f(x,u)=0\qquad\text{in }\Gw,
\end{array}
\ee
satisfying the blow-up condition
\bel{I.2}
\displaystyle
\lim_{d(x)\to 0}u(x)=\infty,
\ee
where  $d(x):=\dist (x,\prt\Gw)$. Whenever a solution to $(\ref{I.1})$-$(\ref{I.2})$ exists it is called a {\it large solution} or an {\it explosive solution}. Although, thanks to \cite[Corollary 3.3]{LGMNA}, in the one-dimensional case $N=1$ with  $f(x,u)\equiv f(u)$ the above problem admits a unique solution, the question of ascertaining whether or not $(\ref{I.1})$-$(\ref{I.2})$ possesses a unique solution received only partial answers  even in the autonomous case when
$f(x,u)\equiv f(u)$ is independent of $x\in\O$. Astonishingly, when $N=1$ the large solution can be unique even when $f(u)$ is somewhere decreasing (see \cite{LGMNA} and \cite{LGMJMAA}), which measures the real level of difficulty of the problem of characterizing the set of $f(x,u)$ for which $(\ref{I.1})$-$(\ref{I.2})$ has a unique positive solution; it is an extremely challenging problem.
\par
Existence of large solutions is associated to the  {\it  Keller--Osserman condition}. When $f$ is independent of $x$, this condition was introduced in \cite{Ke} and \cite{Os} for proving  the first existence results of large solutions in a smooth bounded domain. It reads
\bel{I.3}\begin{array} {lll}
\myint{a}{\infty}\myfrac{ds}{\sqrt{F(s)-F(a)}}<\infty\quad\text{for some }\, a>0\;\text{ where } F(s)=\myint{0}{s}f(t)dt.
\end{array}\ee
When $f=f(x,r)$ a more general version called  in this paper (KO-loc) is introduced in  \cite{LG00} and in  \cite{Ve2}. It asserts that, for any compact subset
 $K$ of $\Gw$, there exists a continuous nondecreasing function $h_K:\BBR_+\mapsto\BBR_+$ such that
 \bel{I.4}\begin{array} {lll}
 f(x,r)\geq h_K(r)\geq 0\quad\text{for all }x\in K\;\text{and }r\geq 0
 \end{array}\ee
 where $h_K$ satisfies
 \bel{I.5}\begin{array} {lll}
\myint{a}{\infty}\myfrac{ds}{\sqrt{H_K(s)-H_K(a)}}<\infty\quad\text{for some }\, a>0\;\text{ where } H_K(s)=\myint{0}{s}h_K(t)dt.
\end{array}\ee
 The condition (KO-loc) guarantees the existence of a  {\it maximal solution, $u^{max}$,} to equation $(\ref{I.1})$. It is obtained as the limit of a decreasing sequence of large solutions $\{u_n\}_{n\in\BBN}$ in an increasing sequence of smooth domains $\{\Gw_n\}_{n\in\N}$ such that $\overline\Gw_n\subset\Gw$ and $\cup_{n\geq 1}\Gw_n=\Gw$ (see e.g. \cite{LG00}, \cite{Ve2}, \cite{MV04}). However, it is not always true that the maximal solution is a large solution. This property depends essentially of the regularity of the domain. If $f(x,u)=u^p$ with $p>1$, the necessary and sufficient condition for such a property to hold is given by \cite{Lab, MV09}. The existence of a {\it minimal large solution} necessitates a minimum of assumptions, either on the regularity of $\Gw$ or on the function $f(x,r)$ (see \cite{Ve2}). Actually, if $\Gw$ is the interior of its closure there exists a decreasing sequence of smooth domains $\Gw'_n$ such that 
$$
  \cap_{n\geq 1} \Gw'_n=\overline\Gw.
$$
If $f$ is defined in $\Gw'\ti \BBR$ where $\Gw'$  is a neighborhood of $\overline \Gw$ with the same monotonicity and (KO-loc) properties therein as in $\overline\Gw\ti\BBR$, and if $f(x,r)\lfloor_{\prt\Gw}>0$, then a sequence of large solutions $\{u_n'\}_{n\in\BBN}$  can be constructed in $\Gw'_n$ and  the limit, $u'$, of the $\{u_n'\}_{n\in\BBN}$ is a candidate for being the minimal large solution, $u^{min}$, since it remains smaller than any large solution in $\Gw$. If $f(x,r)\lfloor_{\prt\Gw}=0$,  then the construction of the minimal solution is possible as soon as for any $n>0$
$(\ref{I.1})$ admits a solution with value $n$ on $\prt\Gw$. For this a minimal regularity condition on $\prt\Gw$ is needed, the Wiener condition \cite[p. 206]{GT}. Furthermore, because of the maximum principle and the fact that $f(x,0)=0$ and $f(x,r)>0$ for all $r>0$ when $x$ belongs to some neighborhood $\CV$ of $\prt\Gw$ the above  (KO-loc) assumption can be weakened in the sense that the function $h_K:\BBR_+\mapsto\BBR_+$ satisfying  $(\ref{I.4})$ and $(\ref{I.5})$ has to exist only when $K$ is a compact subset of $\CV$.
\par
The main property of $u^{max}$ and $u^{min}$ is that any solution $u$ of $(\ref{I.1})$-$(\ref{I.2})$,
 should it exists, satisfies
\begin{equation*}
\begin{array} {lll}
u^{min}(x)\leq u(x) \leq u^{max}(x)\quad\text{for all }\, x\in\Gw.
\end{array}
\end{equation*}
The problem of uniqueness reduces to prove that $u^{max}=u^{min}$. The first results in this direction dealing with $f(x,u)=u^p$  for some  $p>1$, using the asymptotic expansion of any large solution, are proved in \cite{BM}. In this approach, the regularity of the boundary is a crucial assumption. The key point is to prove that
\begin{equation*}
\begin{array} {lll}
\displaystyle\lim_{d(x)\to 0}\myfrac{u^{max}(x)}{u^{min}(x)}=1.
\end{array}
\end{equation*}
After this relation is obtained the uniqueness follows from the fact that there holds
\begin{equation*}
((1+\ge)r)^p\geq  (1+\ge)f(r)\quad\text{for all }\, \ge, r\geq 0.
\end{equation*}
 For regular domains $\Omega$, this technique was substantially refined in \cite{LG06} and  \cite{LG15} to cover the non-autonomous case when $f(x,u)=\mathfrak{a}(x)u^p$ for some non-negative function $\mathfrak{a}(x)$ such that $\mathfrak{a}(x)>0$ for sufficiently small $d(x)$ (see also 
\cite{CRC}, \cite{CRH}, \cite{CD} and \cite{CDG}).
The asymptotic expansion of a large solution near the boundary  requiring so many assumptions, both on the nonlinearity $f$ and on the regularity of $\prt\Gw$, that a new method was introduced in \cite{MV97} in order to bypass this step.  To apply that method the boundary has to satisfy the \emph{local graph condition}, an assumption which is used also in this article. According to it, for every $P\in\prt\Gw$, there exist a neighborhood $Q_P$ of $P$, a positive oriented basis, $\{\vec\nu_1,\ldots,\vec\nu_N\}$, obtained from the canonical one by a rotation, and a function $F\in{C}(\BBR^{N-1};\BBR)$ such that
\begin{equation}
\label{I.6}
\begin{array}{l}
  F(0,\ldots,0)=0,\\[1mm]
  Q_P\cap\O=Q_P\cap\left( \{x\in\BBR^N \;:\; x_N<F(x_1,\ldots,x_{N-1})\}+P\right),
\end{array}
\end{equation}
where the coordinates $(x_1,\ldots,x_n)$ in $(\ref{I.6})$ are expressed with respect to the basis $\{\vec\nu_1,\ldots,\vec\nu_N\}$ (see Figure \ref{Fig1}). Naturally, $\partial \Omega$
satisfies the local graph property if it is Lipschitz continuous.
\begin{figure}[ht]
\centering
\subfigure{\includegraphics[width=10cm]{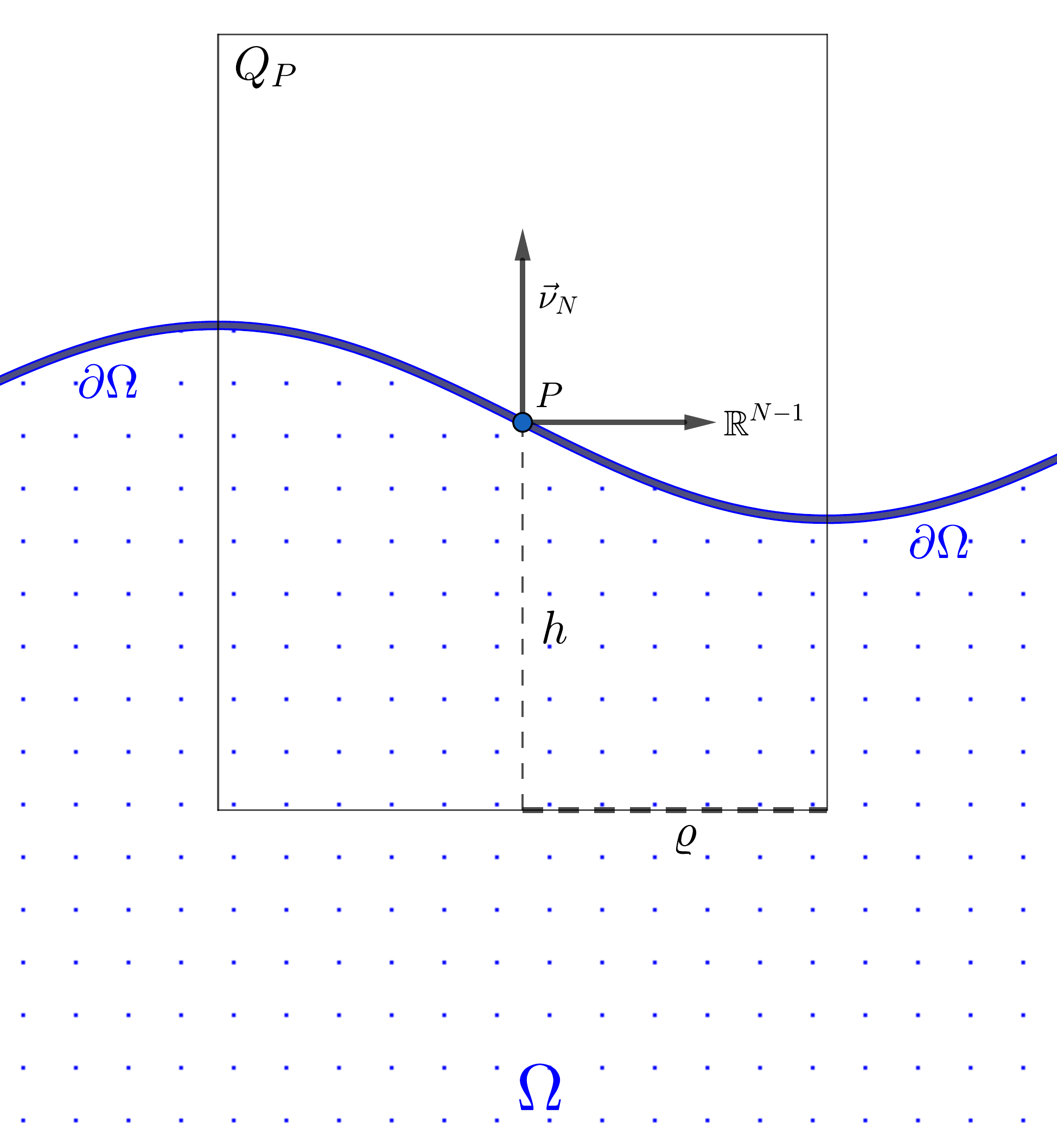}}
\caption{The neighborhood $Q_P$} \label{Fig1}
\end{figure}
\par
 Similarly, in order to avoid the use of the asymptotic expansions  of the large solutions
near the boundary in the proof of the uniqueness, another technique was introduced in \cite{LG07}, and later refined in \cite{CCLG} and \cite{LGMZAMP}, in a radially symmetric context,
based on the strong maximum principle. This technique, which works out even in the context of cooperative systems, \cite{LGM16},  will be combined in this paper
with the technique of \cite{MV97} in order to get the new findings of this paper.
\par
 As far as concerns the nonlinearity $f(\cdot,r)$, in most of the previous papers it is imposed that its rate of decay (or blow-up) near $\prt\Gw$ is a precise function of $d(x)$ (see, e.g., \cite{DH,GLS,LG03,LG06,OX1,OX2,Ve2,Xie09,ZMML10,ZM11}). Throughout this paper it is assumed that  $x\mapsto f(x,r)$ \emph{decays completely nearby} $\prt\Gw$ in the sense that, for every $z\in \prt\Gw$, there exists $\d>0$ such that  $|x-z|<\d$ and $x\in\O$ imply
\begin{equation}\begin{array}{lll}
  \label{I.7}
f(x,\ell+r)-f(x+\ge\vec\nu_N,\ell+r)\geq f(x,\ell)-f(x+\ge\vec\nu_N,\ell)\geq 0\\[1mm]
\hbox{for all }\, r,\ell\geq 0,\;   \hbox{if}\;\; x+\ge \vec\nu_N \in\overline\O\;\;
  \hbox{and}\;\; \ge \in (0,\d).
  \end{array}
\end{equation}
where $\vec\nu_N=(0,...,0,1)$ if $(\ref{I.6})$ holds. Note that this assumption is not  intrinsic to the domain since it depends of the choice of the neighborhood $Q_P$ and the frame $\{\vec\nu_1,\ldots,\vec\nu_N\}$.  In the special case when 
\begin{equation*}
  f(x,r)={\mathfrak a}(x)\tilde f(r)
\end{equation*}
where $\tilde f:\BBR\mapsto\BBR$ is monotone nondecreasing, positive on $(0,\infty)$ and vanishes at $0$, and ${\mathfrak a}\in C(\overline\Gw)$ is nonnegative and positive in a neighborhood of $\prt\Gw$,
the assumption $(\ref{I.7})$ holds if and only if $x\mapsto {\mathfrak a}(x)$ \emph{decays nearby} $\prt\Gw$ in the sense that
\begin{equation}
\label{I.8}
0\leq {\mathfrak a}(x+\ge\vec\nu_N)\leq  {\mathfrak a}(x)\;
 \hbox{ if}\;\; x+\ge \vec\nu_N \in\overline\O\;\;
  \hbox{and}\;\; \ge \in (0,\d).
\end{equation}
\par
If $\Gw$ has a Lipschitz boundary, then there is a truncated circular cone $C_\gg=C\cap B_\gd$ such that any point $P\in\prt\Gw$
is the vertex of the image $\CI_P(C_\gg)$ of $C_\gg$ by an isometry $\CI_P$ of $\BBR^N$ and $\CI_P(C_\gg)\subset\Gw^c$. In such case, $\nu_N$ can be chosen to be the axis of rotational symmetry of $C_\gg$.
\par
 In this paper, associated to $f(x,u)$, we consider the function $g$ defined  on $\overline\Gw\ti\BBR_+$ by
\begin{equation}
\label{I.9}
  g(x,\ell):=\inf\{f(x,\ell+u)-f(x,u) \;:\; u\geq 0\},\quad \text{for all }\,(x,\ell)\in \overline\Gw\ti\BBR_+.
\end{equation}
There always holds $g\leq f$ and $g(x,.)$ is monotone nondecreasing as $f(x,.)$ is. Thus, if $g$ satisfies (KO-loc), so does $f$, but the converse is not true  in general as it is shown in the Appendix.  Moreover, if $f(x,.)$ is convex for all $x\in\overline\Gw$, then $f=g$. This is due to the fact that
\begin{equation*}
\begin{split}
f(x,u+\ell)-f(x,u) & =\myint{0}{\ell} \prt^r_{u} f(x,u+s)ds\\ & \geq\myint{0}{\ell} \prt^r_{u}(x,u'+s)ds=f(x,u'+\ell)-f(x,u')
\end{split}
\end{equation*}
 for all $u\geq u'\geq 0$, $\ell>0$ and $x\in \overline\Gw$, since the right partial derivative $\prt^r_{u}f(x,u)$ of $u\mapsto f(x,u)$  is nondecreasing with $u$. Hence, the minimum of
\[
  u\mapsto f(x,u+\ell)-f(x,u)
\]
is achieved  at $u=0$ and therefore, $f=g$. Finally, if $f$  decays completely nearby $\prt\Gw$, then $g$ also decays in the sense that
 \begin{equation}\begin{array}{lll}
  \label{I.10}
0\leq g(x+\ge\vec\nu_N,r)\leq  g(x,r)\;
\hbox{ for all }\, r\geq 0,\;   \hbox{if}\;\; x+\ge \vec\nu_N \in\overline\O\;\;
  \hbox{and}\;\; \ge \in (0,\d).
  \end{array}
\end{equation}
Furthermore, taking $\ell=0$ in $(\ref{I.7})$, it becomes apparent that $f$  satisfies the same inequality $(\ref{I.10})$ as $g$.\smallskip
\par
The following equation
 \bel{I.11}
-\Gd u+g(x,u)=0\quad\text{in }\,\Gw,
\ee
 closely related to $(\ref{I.1})$ plays a fundamental role in our study.
 Following \cite[Def. 2.6]{MV03}, we introduce the following concept.
\vspace{0.4cm}

\nind{\bf Definition} {\it Let $z\in\prt\Gw$. We say that equation $(\ref{I.11})$
possesses a strong barrier at $z$ if there exists a number $r_z>0$ such that,   for every $r\in (0,r_z]$, there exists a positive supersolution $u=u_{r,z}$ of $(\ref{I.11})$ in $\Gw\cap B_r(z)$   with 
\bel{I.12}\displaystyle
u_{r,z}\in C(\overline\Gw\cap B_r(z))\quad\text{and}\lim_{\tiny\begin{array} {ccc}y\to x\\
y\in \Gw\cap B_r(z)\end{array}}\!\!\!\!\!\!\!\!\!\!\!u_{r,z}(y)=\infty\quad\text{for all }\,x\in \Gw\cap\prt B_r(z).
\ee}
Notice that the local supersolution $u_{r,z}$ of $(\ref{I.11})$
is also a supersolution of $(\ref{I.1})$ since $g\leq f$. Our first result is the following.

\bth{Th1.1}
Suppose that $\Gw$ is Lipschitz continuous and $f\in C({\overline\Gw\ti\BBR})$  satisfies $f(x,0)=0$, $r\mapsto f(x,r)$ is nondecreasing for all $x\in\overline\Gw$, and $f(.,r)$  decays completely nearby $\prt\Gw$ as it is formulated in $(\ref{I.7})$.  Assume, in addition, that  the function $g\in C({\overline\Gw\ti\BBR})$ defined from $f$ by $(\ref{I.9})$ is positive on a neighborhood $\CV$ of $\prt\Gw$ and  satisfies {\rm (KO-loc)}; that is, for any compact subset $K\subset\CV$ there exists a continuous nondecreasing function $h_K:\BBR_+\mapsto\BBR_+$ such that
 \bel{I.13}\begin{array} {lll}
 g(x,r)\geq h_K(r)\geq 0\quad\text{for all }x\in K\;\text{and }r\geq 0,
 \end{array}\ee
where $h_K$ satisfies  $(\ref{I.5})$.  If the equation $(\ref{I.1})$ possesses a strong barrier at any $z\in \prt\Gw$, then the problem $(\ref{I.1})$-$(\ref{I.2})$ possesses a unique solution, i.e.  $ u^{min}=u^{max}$.
\es

The assumption that $g$ satisfies (KO-loc) is actually an assumption on $f$. Indeed, $f$ must grow sufficiently fast at $\infty$ so that $g$ still satisfies (KO-loc). This assumption is weaker than the superadditivity with constant $C$ introduced in \cite{MV04},  according with it 
\bel{I.14}f(x,u+\ell)\geq f(x,u)+f(x,\ell) -C\quad\text{for all }x\in\overline\Gw\, \text{ and }\;u,\ell\geq 0.
\ee
Under the superadditivity assumption,  there holds, for any $\ell,u\geq 0$, that
$$g(x,\ell)\geq f(x,\ell+u)-f(x,u)\geq f(x,\ell)-C.
$$
Therefore, if $f$ satisfies (KO-loc), so does $g$.
Our second  result, valid under a weaker assumption on $\Gw$, requires a new assumption on $f$.

\bth{Th1.2} Assume that $\Gw$ satisfies the local graph property and that the assumptions on $f$ and $g$ in \rth{Th1.1} hold.  Suppose, moreover, that there exists $\phi\in C^2(\BBR_+)$ such that $\phi(0)=0$, $\phi (r)>0$ for $r>0$, and
$$\phi' (r)\geq 0\quad\text{and }\,\phi'' (r)\leq 0\quad \text{for all }\,r\geq 0,$$
for which the function $f$ verifies the inequality
\bel{I.15}
\myfrac{f(x,r+\ge\phi(r))}{f(x,r)}\geq 1+\ge\phi'(r)\quad \text{for all }\,r\geq 0\,\text { and }\; x\in\overline\Gw,
\ee
for some  sufficiently small  $\ge>0$.  Then, the problem $(\ref{I.1})$-$(\ref{I.2})$ possesses at most one solution.
\es

 Although the assumption on $(f,\gf)$ may look unusual, it turns out that  when $\gf (r)=r$  it is equivalent to
\bel{I.16}r\mapsto\myfrac{f(x,r)}{r}\quad \text{is nondecreasing on }\,(0,\infty),
\ee
which is the assumption used in \cite{MV97}. Since  $(\ref{I.16})$ implies that
\begin{equation*}
\myfrac{f(x,r+\ge\ln (1+r))}{r+\ge\ln (1+r)}\geq \myfrac{f(x,r)}{r}
\end{equation*}
or, equivalently,
\bel{I.17}
\myfrac{f(x,r+\ge\ln (1+r))}{f(x,r)}\geq1+\ge\myfrac{\ln (1+r)}{r}\quad \text{for all }\,r\geq 0,
\ee
which entails
\bel{I.18}
\myfrac{f(x,r+\ge\ln (1+r))}{f(x,r)}\geq1+\myfrac{\ge}{1+r}\quad \text{for all }\,r\geq 0,
\ee
 and $(\ref{I.18})$ is $(\ref{I.15})$ for the special choice $\phi(r)=\ln (1+r)$, 
it becomes apparent that  $(\ref{I.15})$ is substantially weaker than $(\ref{I.16})$.

\mysection{Proof of \rth{Th1.1}}

\noindent Since $\O$ is a Lipschitz continuous bounded domain, it satisfies the local graph property at
each point of the boundary. Let $P\in\partial\Omega$ and consider a basis $\{\vec\nu_1,\ldots,\vec\nu_N\}$ and a  neighborhood $Q_P$ satisfying $(\ref{I.6})$. Throughout this proof, it is assumed that any point of $\BBR^N$ is expressed in coordinates with respect to $\{\vec\nu_1,\ldots,\vec\nu_N\}$. Setting
\[
  \hat x:=(x_1,...,x_{N-1})\quad\hbox{for every}\;\;x=(x_1,...,x_{N-1},x_N)\in\BBR^N,
\]
and denoting by $\hat B_\varrho(P)$ the ball of $\BBR^{N-1}$ with center $P=(\hat P,0)$ and radius $\varrho$, we can assume that
\begin{equation*}
  Q_P =\{ x\in\BBR^N\;:\;\; |\hat x -\hat P |<\varrho,\;\; |x_N-P_N|<h\}
  =\hat B_\varrho(P)\ti (P_N-h,P_N+h)
\end{equation*}
for some $\varrho>0, h>0$ such that $\p\O$ is bounded away from the \lq\lq top\rq\rq \, and the \lq\lq bottom\rq\rq\, of $Q_P$ and
\[
 \p\O\cap\overline{Q_P}=\overline{\p\O\cap Q_P}
\]
(see Figure \ref{Fig1}).
\begin{figure}[ht]
\centering
\subfigure{\includegraphics[width=10cm]{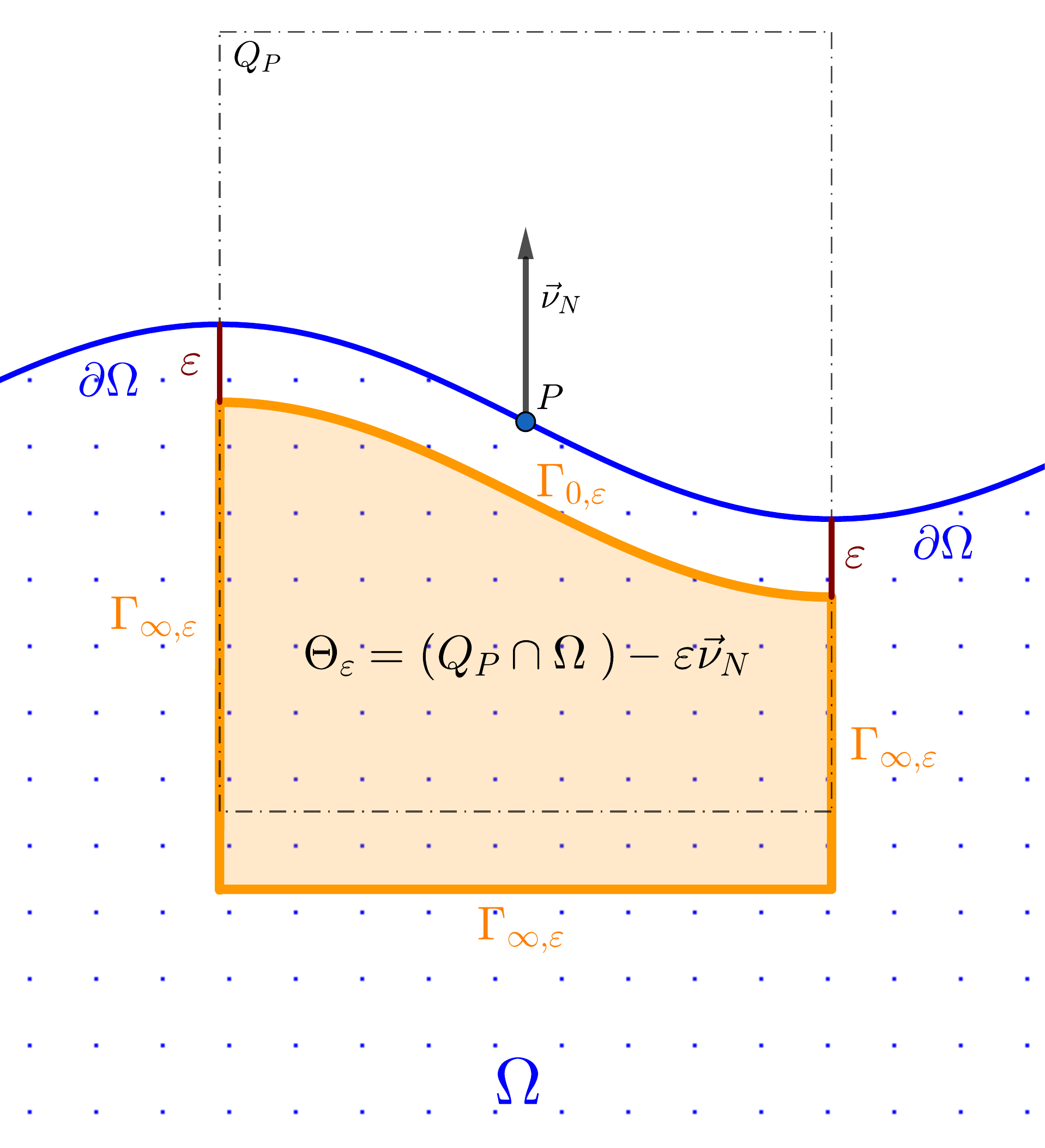}}
\caption{The domain $\T_\ge$} \label{Fig2}
\end{figure}
Thus, setting
\begin{equation}
\label{2.1}
    \T_\ge:=(Q_P\cap \O)-\ge\vec\nu_N,\quad \ge\geq0,
\end{equation}
the existence of  $\ge_1>0$ such that
\begin{equation}
\label{2.2}
  \overline\T_\ge\subset\O\quad\hbox{for all }\;\; 0<\ge<\ge_1
\end{equation}
is guaranteed. Subsequently, we denote
\begin{equation}
\label{2.3}
  \G_{0,\ge}:=(Q_P\cap\p\O)-\ge\vec\nu_N,\qquad \G_{\infty,\ge}:=(\p Q_P\cap \O)-\ge\vec\nu_N,\quad\hbox{for all}\;\;\ge\geq0,
\end{equation}
(see Figure \ref{Fig2}) and consider
\[
  \ge_0:=\min\{\ge_1,\d\},
\]
where $\d$ is the one of $(\ref{I.7})$. Then, the following lemma of technical nature holds.

\blemma{lem2.1} Under the assumptions of \rth{Th1.1}, the problem
\begin{equation}
\label{2.4}
  \left\lbrace\begin{array}{ll}
	-\D \ell +g(x,\ell)=0 \quad&  \text{in }\;\T_0\\
	\phantom{-\D  +g(x,\ell)}
	\ell=0  \quad & \text{in }\;\G_{0,0}\\
	\phantom{-\D  +g(x,\ell)}
	\ell=+\infty \quad& \text{in }\;\G_{\infty,0},
  \end{array}\right.
\end{equation}
admits, at least,  a positive solution, $\ell$.
\es

\noindent \Proof   As we are imposing that $\ell =0$ on $\G_{0,0}$, our singular boundary condition is reminiscent of those
considered previously in \cite{MV97,MV04}.
To construct $\ell$ one can argue as follows. First, consider any increasing sequence of nonnegative  functions, $\{b_n\}_{n\in\N}\subset C^{0,1}(\p \T_{0})$, satisfying
\[
\left\lbrace\begin{array}{ll}\displaystyle
 \phantom{  \lim_{n\rightarrow\infty}} b_n(x)=0 & \quad \hbox{for all}\;\;x\in \G_{0,0}\\\displaystyle
  \lim_{n\rightarrow\infty}b_n(x)=+\infty & \quad \hbox{for all}\;\;x\in \G_{\infty,0},
\end{array}\right.
\]
and let $L_{n}$ be the unique positive solution of the (non-singular) boundary value  problem
\begin{equation}
\label{2.5}
  \left\lbrace\begin{array}{ll}
	-\D L+g(x,L) = 0   \quad&  \text{in }\;\;\T_{0}\\
	\phantom{-\D +g(x,L)}
	L=b_n  \quad&  \text{on }\;\; \p \T_{0}.
  \end{array}\right.
\end{equation}
The solution is the unique minimizer of the lower semicontinuous convex functional
$$J(L)=\myint{\T_{0}}{}\left(\myfrac{1}{2}|\nabla L|^2+G(x,L)\right) dx\quad\text{with }\; G(x,L)=\myint{0}{L}g(x,s) ds,
$$
defined over the affine space of functions in $H^1(\T_{0})$ with trace $b_n$ on $\p \T_{0}$. Since
 $\{b_n\}_{n\in\N}$ is increasing, it follows from the maximum principle that
\[
  L_n\leq L_{n+1} \qquad \hbox{for all}\;\; n\geq 1.
\]
Since $g$ satisfies the strong barrier property there exists $r_{\!\!_P}>0$ such that, for any $r\in (0,r_{\!\!_P}]$, there exists
a supersolution $u_{r,P}$ of $(\ref{I.11})$ in $\Gw\cap B_{r}(P)$ which is continuous in $\overline\Gw\cap B_{r}(P)$. Up to changing $Q_P$ we can assume that for some $r\in (0,r_{\!\!_P}]$,
$$
B_{r}(P)\cap\prt\Gw=Q_P\cap\prt\Gw=\Gg_{0,0}\quad\hbox{and}\quad
\Gw\cap B_{r}(P)\subset Q_P\cap\Gw=\Gth_0.
$$
By the maximum principle, $ L_n\leq u_{r,P}$ in $B_{r}(P)\cap\prt\Gw$. Since $\prt\Gw$ is Lipschitz and the $L_n$ remain locally bounded in a neighborhood of $Q_P\cap\prt\Gw$, it follows by
\cite[Th. 8.29]{GT} that they are locally H\"older continuous near $Q_P\cap\prt\Gw$ and hence the sequence $\{L_n\}_{n\in\N}$ is locally uniformly continuous near $Q_P\cap\prt\Gw$. Therefore, the pointwise limit
\[
 \ell:=\lim_{n\ra\infty}L_{n}
\]
 is well defined in $\Theta_0$ and achieves finite values in $\overline\Gw\cap B_{r}(P)$ since it is dominated by $u_{r,P}$. In what follows we prove that  $\ell$ is continuous in $\overline\Gw\cap Q_P$, vanishes on $\Gg_{0,0}$ and satisfies $(\ref{2.4})$.  For every $\gz\in\T_0$ consider $\tilde{r}>r>0$ so that $\bar{B}_{\tilde{r}}(\gz)\subset\T_0$. Obviously, there exists an integer $n_0$ such that $L_n|_{B_{\tilde{r}}(\gz)}$ is well defined for all $n\geq n_0$. Let  $m$ denote the maximal positive large solution of
\[
  -\D m+g(x,m)=0\quad \hbox{in}\;\; B_{\tilde{r}}(z).
\]
Then, we have that
\begin{equation*}
 0\leq L_n(x)<||m||_{{C}(\ov{B}_r(z))}\quad\hbox{for all}\;\;x\in\bar{B}_r(z)\;\;\hbox{and}\;\; n\geq n_0.
\end{equation*}
Thus, combining a rather standard compactness argument together with the interior Schauder estimates there exists a subsequence,  $\{L_{n_k}\}_{k\in\N}$, which converges locally uniformly to $\ell$ in $\Theta_0$. Clearly $\ell$ satisfies $(\ref{2.4})$, and since the sequence $\{L_n\}_{n\in\N}$ is locally H\"older continuous up to $Q_P\cap\prt\Gw$,  $\ell$ vanishes on $\Gg_{0,0}$.
\qeda
\vspace{0.4cm}

 The next result provides us with a supersolution of $(\ref{I.1})$ in $\Theta_\ge$.

\bprop{pr2.2}

For every $\ge \in (0,\ge_0)$, the function
\begin{equation*}
  \bar{u}_\ge(x)=u^{min}(x+\ge\vec{\nu}_N)+\ell(x+\ge\vec{\nu}_N),\qquad x\in\T_\ge,
\end{equation*}
provides us with a supersolution of $(\ref{I.1})$ in $\T_\ge$ such that
\[
  \bar{u}_\ge=+\infty\quad \hbox{on}\;\; \p\T_\ge.
\]
\es

\noindent\Proof
The fact that $\bar{u}_\ge=+\infty$ on $\p\T_\ge$  follows readily from the definition.  Indeed, by $(\ref{2.1})$ and $(\ref{2.3})$, we have that $x+\ge\vec\nu_N\in\p\O$ if $x\in\p\T_\ge\setminus\G_{\infty,\ge}$. Thus,
\[
  \bar{u}_\ge(x)\geq u^{min}(x+\ge\vec\nu_N)=+\infty\quad\hbox{for all}\;\;x\in\p\T_\ge\setminus\G_{\infty,\ge}.
\]
On the other hand, by $(\ref{2.4})$, we have that, for every $x\in\G_{\infty,\ge}$,
\[
  \bar{u}_\ge(x)\geq\ell(x+\ge\vec{\nu}_N)=+\infty.
\]
Therefore, $\bar{u}_\ge=+\infty$ on $\p\T_\ge$. Now, restricting ourselves to $\T_\ge$, it follows from $(\ref{2.4})$ that
\begin{align*}
  -\D\bar{u}_\ge(x)&=-\D u^{min}(x+\ge\vec{\nu}_N)-\D\ell(x+\ge\vec{\nu}_N)\\
  &=-f(x+\ge\vec{\nu}_N,u^{min}(x+\ge\vec{\nu}_N))-g(x+\ge\vec{\nu}_N,\ell(x+\ge\vec{\nu}_N)).
\end{align*}
Thus, owing to $(\ref{I.7})$-$(\ref{I.10})$,  it becomes apparent that
\begin{align*}
  -\D\bar{u}_\ge(x) &\geq -f(x,u^{min}(x+\ge\vec{\nu}_N))- g(x,\ell(x+\ge\vec{\nu}_N))
\end{align*}
for every $x\in\T_\ge$. Finally, by the definition of $g(x,u)$ (see $(\ref{I.9})$), we find that
\begin{equation*}
\begin{split}
   & -f(x,u^{min}(x\!+\!\ge\vec{\nu}_N))\!-\!g(x,\ell(x\!+\!\ge\vec{\nu}_N))]  \\
   &\geq -[f(x,u^{min}(x\!+\!\ge\vec{\nu}_N))\! +\!f(x,u^{min}(x\!+\!\ge\vec{\nu}_N)\!+\!\ell(x\!+\! \ge\vec{\nu}_N))\!-\! f(x,u^{min}(x\!+\!\ge\vec{\nu}_N))]\\
  &=-f(x,\bar{u}_\ge(x)).
\end{split}
\end{equation*}
Therefore, $\bar{u}_\ge$ is a supersolution of $(\ref{I.1})$ in $\T_\ge$, which ends the proof.
\qeda
\vspace{6mm}

We can complete now the proof of \rth{Th1.1}. By $(\ref{2.2})$, $u^{max}$ is bounded on $\p\T_\ge$ for all $0<\ge\leq \ge_0$. Thus, it follows from the strong maximum principle that
\begin{equation}
\label{2.6}
  \bar{u}_\ge(x)=u^{min}(x+\ge\vec{\nu}_N)+\ell(x+\ge\vec{\nu}_N)\geq u^{max}(x),\quad \hbox{for all}\;\; 0<\ge\leq \ge_0,\;x\in\T_\ge.
\end{equation}
To prove $(\ref{2.6})$ we argue by contradiction. Since
\begin{equation*}
  \bar{u}_\ge(x)=+\infty > u^{max}(x),\quad \hbox{for all}\;\; 0<\ge\leq \ge_0,\;x\in\p\T_\ge,
\end{equation*}
if $(\ref{2.6})$ fails, then, for some $\ge \in (0,\ge_0)$, there exists an open subset, $D=D(\ge)$, with $\overline D \subset \T_\ge$, such that
\begin{equation}
\label{2.7}
  \left\{ \begin{array}{ll}  \bar{u}_\ge=u^{min}(\cdot+\ge\vec{\nu}_N)+\ell(\cdot+\ge\vec{\nu}_N)\lneqq   u^{max}& \quad \hbox{in}\;\; D\\ \bar{u}_\ge = u^{max}& \quad \hbox{on}\;\; \p D.\end{array}\right.
\end{equation}
Thus, setting
\[
  v:= u^{max}-\bar u_\ge,
\]
we find from \rprop{pr2.2} and assumption $(\ref{2.7})$ that
\[
  -\D v = -\D u^{max}+\D \bar u_\ge \leq -[f(x,u^{max})-f(x,\bar u_\ge)]<0 \quad
  \hbox{in}\;\; D,
\]
while $v=0$ on $\p D$. Consequently, $v<0$ in $D$, which implies $u^{max}<\bar u_\ge$ in $D$
and contradicts the assumption $(\ref{2.7})$. This contradiction shows the above claim.
\par
Now, letting $\ge\da 0$ in $(\ref{2.6})$ yields
\begin{equation*}
  u^{min}(x)+\ell(x)\geq u^{max}(x)\quad \hbox{for all}\;\;x\in\T_0.
\end{equation*}
Therefore, it becomes clear that
\begin{equation*}
  0 \geq \limsup_{d(x)\da0}\left( u^{max}(x)-u^{min}(x)\right)\geq 0,
\end{equation*}
which entails
\[
  \lim_{d(x)\da0}\left( u^{max}(x)-u^{min}(x)\right)=0.
\]
Finally, setting
\[
  L:=u^{min}-u^{max}\leq 0,
\]
by the monotonicity of $f$ we find that
\[
  \left\lbrace\begin{array}{ll}
    -\D L=f(x,u^{max})-f(x,u^{min})\geq0 & \hbox{in}\;\; \O\\[1mm]
    \phantom{-\D }
    L=0 & \hbox{on}\;\; \p\O,
  \end{array}\right.
\]
and, consequently, applying the maximum principle,  we can infer that $L=0$.
This ends the proof of \rth{Th1.1}.

\setcounter{equation}{0}
\section{Proof of \rth{Th1.2}}

\noindent We assume that $u^{max}$ is a large solution, i.e. satisfies $(\ref{I.1})$-$(\ref{I.2})$. The next result which has the same expression as  \rlemma{lem2.1} needs actually a slightly different proof due to the
fact that the boundary may not be regular at all.
\blemma{lem3.1} Under the assumptions of \rth{Th1.2}, there exists a nonnegative function
$\ell \in C^1(\Gth_0)$,  bounded on any compact subset of  $\overline\Gw\cap Q_P$, satisfying
\begin{equation}
\label{3.1}
  \left\lbrace\begin{array}{ll}
	-\D \ell +g(x,\ell)=0 \quad&  \text{in }\;\T_0\\
	\phantom{-\D +g(x,\ell)}
	\ell=+\infty \quad& {\color{blue} \text{on }}\;\G_{\infty,0}.
  \end{array}\right.
\end{equation}
\es

\nind\Proof Since the equation $(\ref{I.11})$ admits a strong barrier at $P$,  we can assume that there admits a supersolution in $B_\gr(P)\subset Q_P$, where
$Q_P$ is the cylinder of diameter $\gr$. Hence, $B_\gr(P)\subset Q_P$ and $B_\gr(P)\cap\Gw\subset Q_P\cap\Gw=\Gth_0$. We denote the barrier by
$u_{\gr,P}$. For $\gs>0$ small compared to $\gr$ we consider a domain $\Gth'_{\gs}$ such that
$$
  \Gw\cap \Gth_{\gs}\subset \Gth'_{\gs}\subset \Gw\cap \Gth_{\frac{\gs}{2}},
$$
and we denote by $\Gg'_{0,\gs}$  its upper boundary and by $\Gg'_{\infty,\gs}$  its lateral and lower boundaries. We can assume that $\Gg'_{0,\gs}$ is Lipschitz continuous. Let $\ell=\ell_{n,\gs}$ be the solution, obtained by minimization, of
\begin{equation*}
  \left\lbrace\begin{array}{ll}
	-\D \ell +g(x,\ell)=0 \quad&  \text{in }\;\Gth'_{\gs}\\
	\phantom{-\D +g(x,\ell)}
	\ell=0&  \text{on }\;\Gg'_{0,\gs}
	\\
	\phantom{-\D +g(x,\ell)}
	\ell=n \quad& \text{on }\;\Gg'_{\infty,\gs}.
  \end{array}\right.
\end{equation*}
 Since  $u_{\gr,P}\in C(\overline\Gw\cap B_\gr(P))$ is positive in $\Gw\cap B_\gr(P)$,
 for sufficiently small $\gs$ we have that
$\ell_{n,\gs}\leq u_{\gr,P}$ in $\Gth'_{\gs}\cap B_\gr(P)$.  Thus, By the maximum principle
\begin{equation*}
\ell_{n,\gs}\leq \ell_{n',\gs'}\quad\text{in }\,\Gth'_{\gs}\,\text{ if }\;n'>n\,\text{ and }\,\gs'<\gs.
\end{equation*}
When $\gs\downarrow 0$, $\ell_{n,\gs}$ increases and converges to a function $\ell:=\ell_{n}$ which satisfies
$$  \left\lbrace\begin{array}{ll}
-\D \ell +g(x,\ell)=0 \quad&  \text{in }\;\Gth_{0}
\\
\phantom{-\D +g(x,\ell)}
\ell\leq u_{\gr,P}&  \text{on }\;\Gg_{0,0}
\\
\phantom{-\D +g(x,\ell)}
	\ell=n&  \text{on }\;\Gg_{\infty,0}.
  \end{array}\right.
$$
 As $g$ satisfies (KO-loc), $\ell_n$ remains locally bounded in $\Gth_{0}$.  Therefore, $\ell_{n}\uparrow \ell$ as $n\to\infty$. Clearly,
$\ell$ is bounded on any compact set $K\subset \overline\Gw\cap Q_P$, it belongs to $C^1(\Gth_0)$,  by standard elliptic regularity theory, and satisfies $(\ref{3.1})$.
\qeda
\medskip

\par
Now, suppose that $u(x)$ is any positive solution of $(\ref{I.1})$-$(\ref{I.2})$ and consider
\begin{equation}
\label{3.2}
  \bar u_\e(x):= u(x+\e\vec{\nu}_N)+\ell(x+\e\vec{\nu}_N),\qquad x\in\Theta_\e,
\end{equation}
for sufficiently small $\e>0$. The argument of the proof of \rprop{pr2.2} works out \emph{mutatis mutandis} to show that $\bar u_\e$ is a supersolution of $(\ref{I.1})$  in $\T_\e$. Moreover, by $(\ref{2.2})$, $u$ is bounded on $\p\T_\e$ for sufficiently small $\e>0$. Thus, arguing as in the last step of the proof of \rth{Th1.1}, it follows from the strong maximum principle that
\begin{equation}
\label{3.3}
  \bar{u}_\e(x)=u(x+\e\vec{\nu}_N)+\ell(x+\e\vec{\nu}_N)\geq u^{max}(x),\quad \hbox{for all}\;\; 0<\e\leq \e_0,\;x\in\T_\e.
\end{equation}
As there exists a decreasing sequence $\e_n \da 0$ as  $n\ua +\infty$
such that the function
\[
  \ell=\lim_{n\to \infty}\ell(\cdot +\e_n\vec{\nu}_N)
\]
solves $(\ref{3.1})$, particularizing $(\ref{3.3})$ at $\e=\e_n$ and letting $n\ua +\infty$ yields
\begin{equation}
\label{3.4}
  u(x)+\ell(x)\geq u^{max}(x)\quad \hbox{for all}\;\;x\in\T_0.
\end{equation}
On the other hand, by the definition of $u^{max}$ there holds
\begin{equation*}
  u^{max}(x)+\ell(x)\geq u(x)\quad \hbox{for all}\;\;x\in\T_0.
\end{equation*}
Therefore, for every $x\in\T_0$, we have that
\begin{equation}
\label{3.5}
  \ell(x) \geq u^{max}(x)-u(x) \geq 0.
\end{equation}
Finally, in order to infer from  $(\ref{3.5})$ that $u(x)=u^{max}(x)$ for all $x\in\T_0$, we will use the next result of technical nature.

\blemma{lem3.2}
Let $u_1(x)$ and $u_2(x)$ be positive solutions of $(\ref{I.1})$-$(\ref{I.2})$ such that
\begin{equation}
\label{3.6}
\lim_{d(x)\da 0}\frac{u_2(x)-u_1(x)}{\v(u_1(x))}=0.
\end{equation}
Then, $u_1=u_2$ in $\O$.
\es
\Proof
For sufficiently small $\e>0$, consider the function $v$ defined by
\begin{equation}
\label{3.7}
  v:=u_1+\e \v(u_1),
\end{equation}
where $\v$ is the function introduced in the statement of \rth{Th1.2}.
We claim that $v\geq u_2$ in a neighborhood of $\p\O$. Indeed, by $(\ref{3.6})$, for any $\ge>0$ there exists $\gd>0$ such that if $d(x)<\d$, then
\[
   \frac{u_2(x)-u_1(x)}{\v(u_1(x))} \leq \e = \frac{v(x)-u_1(x)}{\v(u_1(x))}.
\]
Thus, $v(x)\geq u_2(x)$ provided $d(x)\leq \d$. On the other hand, since $\v''\leq 0$, we have that
\begin{align*}
   -\D v & = -\D u_1 -\e\v'(u_1)\D u_1-\e\v''(u_1) |\nabla u_1|^2 \\[5pt] & \geq
   -(1+\e\v'(u_1))\D u_1 \\[5pt] & = -(1+\e\v'(u_1))f(x,u_1(x)).
\end{align*}
Hence,
\[
  -\D v + f(x,v) \geq f(x,v)-(1+\e\v'(u_1))f(x,u_1).
\]
Consequently, thanks to $(\ref{3.7})$ and $(\ref{I.15})$,  it is clear that 
\[
  -\D v + f(x,v) \geq 0
\]
in $\O$. So, $v$ is a supersolution of $(\ref{I.1})$ and hence, $v\geq u_2$ in $\O$ for sufficiently small $\e>0$. Thus, letting $\e\da 0$ yields $u_1\geq u_2$ in $\O$. By symmetry, $u_1=u_2$ holds,
 which ends the proof.
\qeda
\vspace{0.4cm}

Dividing $(\ref{3.5})$ by $\v(u(x))$ and letting $d(x)\da 0$, yields
\[
\lim_{d(x)\da 0}\frac{u^{max}(x)-u(x)}{\v(u(x))}=0.
\]
Consequently, by \rlemma{lem3.2}, we find that $u=u^{max}$. This ends the proof of \rth{Th1.2}.

\setcounter{equation}{0}
\section{Appendix}
\subsection{On the Keller-Osserman condition}

 The next result shows how imposing the Keller--Osserman condition on the associated function $g$ is stronger than imposing it on $f$.

\bprop{xx&}  There are increasing functions $f$ that satisfy  (KO) and such that the corresponding function $g$ does not. \es
\Proof To construct such an example, one can consider any function $f$ such that
\[
  u^2\leq f(u)\leq u^3 \quad \hbox{and}\quad f(u)=f(\min I_n) \quad \hbox{for all}\;\; u\in I_n,
\]
where $I_n$, $n\geq 1$, is an arbitrary sequence of intervals such that
\[
  \lim_{n\to +\infty} (\max I_n-\min I_n)  =+\infty\quad \hbox{and}\quad
  \max I_n <\min I_{n+1}\quad \hbox{for all}\;\; n\in\N.
\]
By the properties of $u^2$ and $u^3$, such  a sequence of intervals exists. For this choice we have that, for any given $\ell>0$ and $u>0$, $[u,\ell+u]\subset I_n$ for sufficiently large $n>0$ and hence,
\[
  f(\ell+u)-f(u)=0.
\]
Thus, $g(\ell)=0$. Therefore, $g\equiv 0$, which does not satisfy (KO).

\subsection{On the strong barrier property}

The general problem of finding  conditions so that the strong barrier property occurs  is open. We give below some cases where it holds and a case where it does not. They all deal with nonlinearity of the form
\bel{4.1}
f(x,r)={\mathfrak a}(x)\tilde f(r)
\ee
where ${\mathfrak a}\in C(\overline\Gw)$ is nonnegative and positive near $\prt\Gw$ and $\tilde f:\BBR_+\mapsto\BBR_+$ is continuous and nondecreasing, vanishes at $0$ and satisfies $(\ref{I.3})$.
\medskip

\nind 1- If ${\mathfrak a}>0$ on $\prt\Gw$, then the Keller--Osserman condition holds in $\overline\CV$, where $\CV$ is a neighborhood of
$\prt\Gw$,  because  the function ${\mathfrak a}$ can be extended to $\Gw^c$ as a continuous and positive function by Whitney embedding theorem (see e.g. \cite{EG}). It is a completely open problem to find {\color{blue} out} sufficient conditions in the case where ${\mathfrak a}>0$ vanishes on the boundary.\smallskip
\medskip

\nind 2- If $\prt\Gw$ is $C^2$ and,  for some $\ga>0$,
$$
  g(x,r)\geq d^\ga(x)u^p
$$
it is proved in \cite{MV03} that the strong barrier property holds.
When $\prt\Gw$ is Liptschiz  the distance function loses its  intrinsic interest and has often to be replaced by the first eigenfunction $\phi_1$ of $-\Gd$ in $H^{1}_0(\Gw)$.  In such case, we conjecture that the strong barrier property holds if
$$
  g(x,r)\geq \phi_1^\ga(x)u^p
$$
for some $\a>0$.
\medskip

\nind 3-  If $\prt\Gw$ is $C^2$ and
$$
   g(x,r)\leq e^{-\frac{\gk}{d(x)}}r^p
$$
with $\gk>0$ and $p>1$, then the strong barrier property does not hold. Indeed, it is proved in \cite{MV04-1} that, for every $a\in\partial\Omega$ and $k>0$,  the problem
\begin{equation}
\label{4.2}\begin{array} {lll}
-\Gd u+e^{-\frac{\gk}{d(x)}}u^p=0\quad&\text{in }\;\Gw,\\[0mm]
\phantom{-\Gd +e^{-\frac{\kappa}{d(x)}}u^p}
u=k\gd_a\quad& \text{on }\;\prt\Gw,
\end{array}
\end{equation}
 admits a unique positive solution, $v_{a,k}$. Furthermore, the nonlinearity $r\mapsto r^p$ satisfies the Keller--Osserman condition.  Hence, the equation
\bel{4.3}-\Gd u+e^{-\frac{\gk}{d(x)}}u^p=0\quad\text{in }\;\Gw
\ee
admits a minimal, $u^{min}$, and a maximal, $u^{max}$, large solution (probably they are equal). However,  $v_{a,k}\uparrow u^{min}$ when $k\to\infty$.  Arguing by contradiction, assume that the equation satisfies the strong barrier property at $z\in\prt\Gw$. Then, there exists  $r>0$ such that
the solution $u:=u_n$ of the problem
$$\begin{array} {lll}
-\Gd u+e^{-\frac{\gk}{d(x)}}u^p=0\quad&\text{in }\;B_r(z)\cap \Gw,
\\[0mm]
\phantom{-\Gd +e^{-\frac{\gk}{d(x)}}u^p}
u=n\quad&\text{on }\;\Gw\cap \prt B_r(z),\\[0mm]
\phantom{-\Gd +e^{-\frac{\gk}{d(x)}}u^p}
u=0\quad&\text{on }\;\prt\Gw\cap  B_r(z),
\end{array}
$$
converges, as $n\to\infty$, to a barrier function $u_{r,z}\in C(\overline\Gw\cap  B_r(z))$ satisfying

$$\begin{array}{lll}
-\Gd u+e^{-\frac{\gk}{d(x)}}u^p=0\quad&\text{in }\;B_r(z)\cap \Gw,
\\[0mm]
\phantom{-\Gd +e^{-\frac{\gk}{d(x)}}u^p}
u=\infty\quad&\text{on }\;\Gw\cap \prt B_r(z).
\end{array}
$$
Taking a point $a\in \prt\Gw\cap B^c_{2r}(z)$, for any $k>0$ there exists $n=n(k)$ such that $v_{a,k}\leq n(k)$ on $\Gw\cap \prt B_r(z)$. Since
$v_{a,k}=0$ on $\prt\Gw\cap  B_r(z)$, it follows that $v_{a,k}\leq u_n$.   Thus, letting $k\to\infty$, yields  $u^{min}\leq u_{r,z}$, which is a contradiction. \smallskip
\medskip

\nind 4-  If $\prt\Gw$ is $C^2$ and
$$
  g(x,r)= e^{-\frac{1}{d^\ga(x)}}r^p,
$$ 
with $0<\ga<1$ and $p>1$, it is proved in \cite{SV1} that the limit when $k\to\infty$ of the solutions $v_{a,k}$ of
\begin{equation}
\label{4.4}\begin{array} {lll}
-\Gd u+e^{-\frac{1}{d^\ga(x)}}u^p=0\quad&\text{in }\;\Gw\\[0mm]
\phantom{-\Gd +e^{-\frac{1}{d^\ga(x)}}u^p}
u=k\gd_a\quad&\text{on }\;\prt\Gw,
\end{array}
\end{equation}
is a solution of
$$-\Gd u+e^{-\frac{1}{d^\ga(x)}}u^p=0\quad\text{in }\;\Gw
$$
which vanishes on $\prt\Gw\setminus\{a\}$ and blows up at $a$.
We conjecture that the strong barrier property holds if
$$
  g(x,r)\geq  e^{-\frac{1}{d^\ga(x)}}r^p.
$$


\end{document}